\documentclass[twosided]{degm}
\usepackage{amssymb}
\usepackage{amsmath}
\usepackage{latexsym}
\usepackage{amsthm}
\usepackage{graphicx,pst-node,pst-plot,pstricks}
\newcommand{\LastUpdate}{14 January 2002}

\newcommand{\ol}[1]{\overline{#1}}
\newcommand{\bom}[1]{\boldsymbol{#1}}
\newcommand{\card}{\mathop{\rm card}}

\newcommand{\diag}{\mathop{\rm diag}}

\newcommand{\id}{\mbox{\rm1\hspace{-.2ex}\rule{.1ex}{1.44ex}}
   \hspace{-.82ex}\rule[-.01ex]{1.07ex}{.1ex}\hspace{.2ex}}

\newcommand{\elaw}{{\stackrel{\textrm{d}}{=}}}

\newcommand{\Rk}{{\par\noindent{\bf Remark:~}}}

\newcommand{\saw}{{\textrm{saw}}}

\newcommand{\cB}{{\cal B}}

\newcommand{\cF}{{\cal F}}
\newcommand{\cG}{{\cal G}}
\newcommand{\cH}{{\cal H}}

\newcommand{\cK}{{\cal K}}
\newcommand{\cL}{{\cal L}}

\newcommand{\cO}{{\cal O}}

\newcommand{\bM}{{\bf M}}

\newcommand{\bu}{{\bf u}} 
\newcommand{\bv}{{\bf v}}

\newcommand{\BbA}{\mathbb{A}}

\newcommand{\BbE}{\mathbb{E}}

\newcommand{\BbG}{\mathbb{G}}
\newcommand{\BbH}{\mathbb{H}}

\newcommand{\BbL}{\mathbb{L}}

\newcommand{\BbN}{\mathbb{N}}
\newcommand{\BbO}{\mathbb{O}}
\newcommand{\BbP}{\mathbb{P}}

\newcommand{\BbR}{\mathbb{R}}
\newcommand{\BbS}{\mathbb{S}}

\newcommand{\BbV}{\mathbb{V}}

\newcommand{\BbX}{\mathbb{X}}

\newcommand{\BbZ}{\mathbb{Z}}

\makeatletter
\def\thmhead@plain#1#2#3{%
  \thmname{#1}\thmnumber{\@ifnotempty{#1}{ }#2}%
  \thmnote{ \the\thm@notefont(#3)}}
\let\thmhead\thmhead@plain
\def\swappedhead#1#2#3{%
  \thmnumber{#2}\thmname{\@ifnotempty{#2}{. }#1}%
  \thmnote{ \the\thm@notefont(#3)}}
\makeatother

\theoremstyle{definition} 

 \newtheorem{definition}{Definition}[section]
 \newtheorem{remark}[definition]{Remark}
 \newtheorem{example}[definition]{Example}

\theoremstyle{plain}      

 \newtheorem{theorem}[definition]{Theorem}
 
 \newtheorem{lemma}[definition]{Lemma}

\def\card{{\mathop{\rm card}}}

\begin{document}

\title{\protect{\textsf{On the physical relevance of random walks: 
an example of random walks on a randomly oriented lattice}}
\footnote{Extended version of the talk given by at the Workshop ``Random walks
and geometry'' held at Erwin Schr\"odnger Insitut, Vienna 18 June -- 13 July 2001.\\
\textit{1991 Mathematics Subject Classification:}
60J10, 60K20\\
\textit{Key words and phrases:}
Markov chain, random environment, recurrence criteria, random graphs, oriented graphs, $C^*$-algebras, quantum
communication.
 } }

\author{M.\ Campanino, D.\ Petritis}

\markboth{M.\ Campanino, D.\ Petritis}{Random walks on randomly oriented lattices}

\maketitle

\centerline{ \LastUpdate}
\vskip4mm

\begin{abstract}
Random walks on general graphs play an important role in 
the understanding of the general theory of stochastic processes.
Beyond their fundamental interest in probability theory, they arise
also as simple models of physical systems. 
A brief survey of the physical relevance of the notion of random walk on both
undirected and directed graphs is given followed by the exposition
of some recent results on random walks on randomly oriented lattices.
 It is worth noticing that general
undirected graphs are associated with (not necessarily Abelian) groups
while directed graphs are associated with (not necessarily
Abelian) $C^*$-algebras. Since quantum mechanics is naturally formulated
in terms of  $C^*$-algebras, the study of random walks on directed lattices
has been motivated lately by the development of the new field of quantum
information and communication. 
\end{abstract}

\section{Introduction}
Random walks are mathematical objects thoroughly studied nowadays by probabilists
for their own interest but also for shedding new light into a variety of mathematical
problems like diffusions on manifolds, harmonic analysis, infinite graph theory, group
theory, etc. The other contributions in this volume give an overview of the
most recent developments of random walks in connection with most of these
mathematical disciplines.

An interesting class of problems concerns random walks in random environments. The present contribution
intends to present some recent results on a particular kind of random environment defined by
the orientation of some edges of the graph on which the random walk evolves. However, we felt that
instead of merely announcing these results and reproducing the proofs --- that can anyway be found 
elsewhere \cite{CamPet-rwrol} --- it should be useful for the mathematical community to enlarge this report by
giving
a short overview on the physical relevance of random walks
both on undirected lattices  and
directed lattices (based on some recent connections between them and quantum evolution,
established recently in
\cite{Ler} and in \cite{Pet-QI}).  Therefore, this contribution
is organised as follows. In section \ref{sec-notations} we give the notation, definitions, and
our main results on the probabilistic problem we have studied \cite{CamPet-rwrol}, 
in section \ref{sec-rwul} we give
a short overview of the usefulness of random walks on undirected lattices as representations
of physical quantities like Green's functions, in section \ref{sec-rwol} we 
recall an algebraic construction of (partially) oriented
lattices in terms of $C^*$-algebras and demonstrate the usefulness of random walks in a particular
example stemming from quantum communication. Finally in section \ref{sec-proofs}
we sketch the proof of our results; for detailed proofs the reader can consult \cite{CamPet-rwrol}.

\section{Notation, definitions, and main results}
\label{sec-notations}
We denote $\mathbb{N}=\{0,1,2,\ldots\}$ and $\mathbb{N}_+=\{1,2,3,\ldots\}$. An 
\textit{oriented} (or equivalently \textit{directed}) graph $\mathbb{G}$ is the
quadruple $\mathbb{G}=(G^0,G^1,r,s)$, where $G^0$ and $G^1$ are
denumerable sets called respectively \textit{vertex} and \textit{edge} sets
and $r,s$ are mappings $G^1\rightarrow G^0$ called respectively \textit{range} and
\textit{source} functions. For every $v\in G^0$, the set $N_v^-=s^{-1}(\{v\})\subseteq G^1$ represents
the set of outgoing edges from vertex $v$ and its cardinality $d_v^-=\card N_v^-$ is called
the \textit{outwards degree} of $v$. We assume that our graphs are \textit{row-finite} i.e.\ 
$d_v^-<\infty, \forall v\in G^0$. Similarly, we denote
$N_v^+=r^{-1}(\{v\})\subseteq G^1$ the set of incoming edges at vertex $v$ and $d_v^+=\card N_v^+$
the \textit{inwards degree} of $v$. We assume that our graphs are \textit{locally finite}
i.e.\ $d_v^+<\infty, \forall v\in G^0$. Finally, we denote $N_v=N_v^+\cup N_v^-$ and $d_v=\card N_v$.
The sets $G^0$ and $G^1$ are primitive objects for the graph $\mathbb{G}$. For every
$n\in \mathbb{N}_+$, we define a sequence of derived objects 
\[G^n=\{\alpha=(a_1,\ldots, a_n): a_i\in G^1, i=1,\ldots,n \ \textrm{and}\ 
r(a_i)=s(a_{i+1}), i=1,\ldots, n-1\},\]
called \textit{oriented paths} of length $n$, i.e.\ composed of $n$ edges of $G^1$.
For $\alpha\in G^n$, we denote $|\alpha|=n$ its length. The mappings $r,s$, initially
defined on $G^1$ are naturally extended to $G^n$:  
for $\alpha=a_1\ldots a_n$, $r(\alpha)\equiv r(a_n)$ and $s(\alpha)\equiv s(a_1)$.
The vertices $v\in G^0$ are considered as paths of length  0 and if $\alpha\in G^0$
then $r(\alpha)=s(\alpha)=\alpha$.

The set $G^*=\cup_{n\in\mathbb{N}} G^n$ represents the set of oriented paths of
arbitrary (finite) lenght and 
\[\partial G^*\equiv G^\infty=\{\alpha= (a_i)_{i=0}^\infty: s(a_{i+1})=r(a_i), i\in \mathbb{N}\}\]
represents the set of infinite paths. Notice that thus the set of oriented paths of $\mathbb{G}$
acquires a natural tree structure.

All graphs considered here are assumed \textit{transitive}, i.e.\ for
every $u,v\in G^0$, there is an $\alpha \in G^*$ such that $s(\alpha)=u$ and $r(\alpha)=v$,
Transitivity implies in particular the \textit{no-sink} condition $d_v^-\geq 1$, for all $v\in G^0$.

\begin{remark}
Notice that the above definition of $\mathbb{G}$ is quite general.
\begin{itemize}
\item
It does not exclude multiple edges since no conditions are imposed on the mappings
$r,s$.  If both functions are injective then all edges are simple and moreover each edge
can be thought as a pair of vertices, i.e.\ in that case $G^1\subseteq G^0\times G^0$.
\item
It does not exclude loops (of length 1) since it is not excluded that $r(\alpha)=s(\alpha)$ for
some $\alpha \in G^1$. We say that the graph has no loops if $r(\alpha)\ne s(\alpha)$ for all
$\alpha\in G^1$.
\item 
Finally, the existence of unoriented edges is not excluded since an unoriented edge
can be thought as a pair of oriented edges $\alpha, \beta \in G^1$ with $r(\alpha)=s(\beta)$ and
$r(\beta)=s(\alpha)$.
\end{itemize}
\end{remark}

\begin{definition}
Let $\mathbb{G}=(G^0,G^1,r,s)$ be an oriented, transitive, row-finite graph. A \textit{simple
random walk} on $\mathbb{G}$ is a $G^0$-valued Markov chain 
$(\bM_n)_{n\in\mathbb{N}}$ with stochastic
matrix defined by
\[P_{u,v}=\mathbb{P}(\bM_{n+1}=v|\bM_n=u)=\left\{\begin{array}{ll}
\frac{1}{d_u^-} & \textrm{if} \  \exists \alpha \in G^1: r(\alpha)= v, s(\alpha)=u\\
0 & \textrm{otherwise.}
\end{array}\right.\]
\end{definition}

In the sequel, the graphs we shall consider  are always transitive, row-finite, without multiple
edges and without loops (of length 1). To simplify notation, 
we denote $\mathbb{V}\equiv G^0$ and
$\mathbb{A}\equiv G^1\subseteq \mathbb{V}\times \mathbb{V} \setminus\diag(\mathbb{V}\times \mathbb{V})$.
The range and source functions for the edges are then 
naturally determined and we denote the graph simply 
$\mathbb{G}=(\mathbb{V},\mathbb{A})$.

\begin{definition}
Let $\mathbb{V}=\mathbb{V}_1\times\mathbb{V}_2=\mathbb{Z}^2$, 
with $\mathbb{V}_1$ and $\mathbb{V}_2$ isomorphic to $\mathbb{Z}$
and $\bom{\epsilon}=(\epsilon_y)_{y\in\mathbb{V}_2}$ be a sequence of
$\{-1,1\}$-valued variables assigned to each
ordinate. We call $\bom{\epsilon}$-\textit{horizontally oriented lattice}
$\BbG=\BbG(\BbV,\bom{\epsilon})$,
the directed graph with vertex set $\BbV=\BbZ^2$ and edge set $\BbA$
defined by the condition
$(\bu,\bv)\in\BbA$ if, and only if, $\bu$ and $\bv$ are distinct vertices
satisfying one of the following conditions:
\begin{enumerate}
\item
either $v_1=u_1$ and $v_2=u_2\pm1$,
\item
or $v_2=u_2$ and $v_1=u_1+\epsilon_{u_2}$.
\end{enumerate}
\end{definition}

\begin{example}[The Alternate lattice $\BbL$]
In that case, $\bom{\epsilon}$ is the deterministic sequence $\epsilon_y=(-1)^y$ for
$y \in\BbV_2$.
The figure \ref{fig-alternate} depicts a part of this graph.
\end{example}

\begin{figure}[h]
\centerline{%
\hbox{%
\psset{unit=4.5mm}
\pspicture(-0.5,-0.5)(7,7)
\Cartesian(4.5mm,4.5mm)
\psline[linewidth=0.1pt]{<-<}(-0.5,0)(6.5,0)
\psline[linewidth=0.1pt]{>->}(-0.5,1)(6.5,1)
\psline[linewidth=0.1pt]{<-<}(-0.5,2)(6.5,2)
\psline[linewidth=0.1pt]{>->}(-0.5,3)(6.5,3)
\psline[linewidth=0.1pt]{<-<}(-0.5,4)(6.5,4)
\psline[linewidth=0.1pt]{>->}(-0.5,5)(6.5,5)
\psline[linewidth=0.1pt]{<-<}(-0.5,6)(6.5,6)
\psline[linewidth=0.1pt]{-}(0,-0.5)(0,6.5)
\psline[linewidth=0.1pt]{-}(1,-0.5)(1,6.5)
\psline[linewidth=0.1pt]{-}(2,-0.5)(2,6.5)
\psline[linewidth=0.1pt]{-}(3,-0.5)(3,6.5)
\psline[linewidth=0.1pt]{-}(4,-0.5)(4,6.5)
\psline[linewidth=0.1pt]{-}(5,-0.5)(5,6.5)
\psline[linewidth=0.1pt]{-}(6,-0.5)(6,6.5)
\rput[bl](2.5,3.2){\tiny{$(0,0)$}}
\endpspicture}}
\caption{\label{fig-alternate}
\textsf{The alternately directed lattice $\BbL$ corresponding to the choice
 $\epsilon_{y}=(-1)^{y}$.}}
\end{figure}

\begin{example}[The half-plane one-way lattice $\BbH$]
Here $\bom{\epsilon}$ is the deterministic sequence
\[\epsilon_y=\left\{\begin{array}{ll}
1 & \textrm{if } \ \ y\geq 0\\
-1 & \textrm{if } \ \ y<0.
\end{array}
\right.\]
The figure \ref{fig-halfplane} depicts a part of this graph.
\end{example}
\begin{figure}[h]
\centerline{%
\hbox{%
\psset{unit=4.5mm}
\pspicture(-0.5,-0.5)(7,7)
\Cartesian(4.5mm,4.5mm)
\psline[linewidth=0.1pt]{<-<}(-0.5,0)(6.5,0)
\psline[linewidth=0.1pt]{<-<}(-0.5,1)(6.5,1)
\psline[linewidth=0.1pt]{<-<}(-0.5,2)(6.5,2)
\psline[linewidth=0.1pt]{>->}(-0.5,3)(6.5,3)
\psline[linewidth=0.1pt]{>->}(-0.5,4)(6.5,4)
\psline[linewidth=0.1pt]{>->}(-0.5,5)(6.5,5)
\psline[linewidth=0.1pt]{>->}(-0.5,6)(6.5,6)
\psline[linewidth=0.1pt]{-}(0,-0.5)(0,6.5)
\psline[linewidth=0.1pt]{-}(1,-0.5)(1,6.5)
\psline[linewidth=0.1pt]{-}(2,-0.5)(2,6.5)
\psline[linewidth=0.1pt]{-}(3,-0.5)(3,6.5)
\psline[linewidth=0.1pt]{-}(4,-0.5)(4,6.5)
\psline[linewidth=0.1pt]{-}(5,-0.5)(5,6.5)
\psline[linewidth=0.1pt]{-}(6,-0.5)(6,6.5)
\rput[bl](2.5,3.2){\tiny{$(0,0)$}}
\endpspicture}}
\caption{\label{fig-halfplane}
\textsf{The half-plane one-way  lattice $\BbH$ with
$\epsilon_y=-1$, if $y<0$ and
$\epsilon_y=1$, if $y\geq0$.}}
\end{figure}
\begin{example}[The lattice with random horizontal
orientations $\BbO_{\bom{\epsilon}}$]
Here $\bom{\epsilon}=(\epsilon_y)_{y\in\BbV_2}$
is a sequence of Rademacher, \textit{i.e.}\
 $\{-1,1\}$-valued
symmetric Bernoulli random variables,
that are independent for different values of $y$.
The figure \ref{fig-randomly} depicts part of a realisation of this graph.
The random sequence  $\bom{\epsilon}$ is also termed the
\textit{environment of random horizontal
directions}.
\end{example}

\begin{figure}[h]
\centerline{%
\hbox{%
\psset{unit=4.5mm}
\pspicture(-0.5,-0.5)(7,7)
\Cartesian(4.5mm,4.5mm)
\psline[linewidth=0.1pt]{<-<}(-0.5,0)(6.5,0)
\psline[linewidth=0.1pt]{<-<}(-0.5,1)(6.5,1)  
\psline[linewidth=0.1pt]{>->}(-0.5,2)(6.5,2)
\psline[linewidth=0.1pt]{<-<}(-0.5,3)(6.5,3)
\psline[linewidth=0.1pt]{>->}(-0.5,4)(6.5,4)
\psline[linewidth=0.1pt]{>->}(-0.5,5)(6.5,5)
\psline[linewidth=0.1pt]{<-<}(-0.5,6)(6.5,6)
\psline[linewidth=0.1pt]{-}(0,-0.5)(0,6.5)
\psline[linewidth=0.1pt]{-}(1,-0.5)(1,6.5)
\psline[linewidth=0.1pt]{-}(2,-0.5)(2,6.5)
\psline[linewidth=0.1pt]{-}(3,-0.5)(3,6.5)
\psline[linewidth=0.1pt]{-}(4,-0.5)(4,6.5)
\psline[linewidth=0.1pt]{-}(5,-0.5)(5,6.5)
\psline[linewidth=0.1pt]{-}(6,-0.5)(6,6.5) 
\rput[bl](2.5,3.2){\tiny{$(0,0)$}}
\endpspicture}}
\caption{\label{fig-randomly}
\textsf{The randomly horizontally
directed  lattice $\BbO_{\bom{\epsilon}}$ with
$(\epsilon_y)_{y\in\BbZ}$ an independent and identically distributed
sequence of Rademacher random variables.}}  
\end{figure}

Now we can state our main results.

\begin{theorem}
\label{th-L}
The simple random walk on the alternate lattice $\BbL$ is recurrent.
\end{theorem}

\begin{remark} This result can be easily generalised to any lattice
with  periodically alternating  horizontal directions (for every
finite period).
\end{remark}

\begin{theorem}
\label{th-H}
The simple random walk on the half-plane one-way lattice $\BbH$ is transient.
\end{theorem}

\begin{remark}
 The result concerning transience in theorem \ref{th-H} is robust.
In particular, perturbing the orientation of any finite set of
horizontal lines either by reversing the orientation of these
lines or by transforming them into two-ways does not change the
transient behaviour of the simple random walk.
Therefore, the half-plane one-way lattice is so deeply in the
transience region that the asymptotic behaviour of the
simple random walk cannot be changed by
simply modifying the transition probabilities along
a lower dimensional manifold as  was
the case in \cite{MenPet-wedge} where the bulk
behaviour is on the critical point and it can
be changed by lower-dimensional perturbations.
\end{remark}

\begin{theorem}
\label{th-O}
For almost all realisations of the environment $\bom{\epsilon}$,
the simple random
walk on the randomly horizontally oriented lattice
$\BbO_{\bom{\epsilon}}$ is transient and its speed is $0$.
\end{theorem}

\section
{The physical relevance of random walks on undirected lattices}
\label{sec-rwul}
\subsection
{A brief history}
\label{ssec-history}
The original impetus for the study of the continuous time analogue
of a random walk, the Brownian motion, was given by the seminal work of
Einstein \cite{Ein} on diffusions\footnote{A more easily accessible 
source than the original paper \cite{Ein}, exposing the main ideas
in an informal but fascinating style, is the commented scientific
biography of Einstein \cite{Pai}.}.

The discrete time process, we nowadays call simple random walk,
was first studied by P\'olya \cite{Pol}.
It is remarkable however that, contrary to Brownian motion whose physical
motivation lies in the very definition of the model, the intrinsic physical
relevance of random walks was discovered much later, with the development
of polymer physics \cite{Flo}. Polymers are
long, \textit{topologically} one dimensional molecules,
composed by repeating several times (typically 100--10000) the same
structural unit. The structural unit can be viewed as a small straight
and rigid segment that can be glued with subsequent segments by loose
bonds in such a way that if the first segment is held fixed, the second
segment forms with the previous one a given angle $\theta$ (depending
only on the chemical nature of the molecule) but otherwise its position
is arbitrary. Assuming that the structural units have unit length, they
merely represent directions $\hat{x}\in\BbS^2$. Hence, denoting
$C(\hat{x}, e)$ the cone of opening $2\theta$, having its apex at the end
point, $e$, of a segment and its axis collinear with $\hat{x}$, the second
segment  will lie on an arbitrary separatrix of $C(\hat{x},e)$. One thus
immediately recognises a $\BbR^3$-valued discrete time random process
$ (S_n)_{n\in\BbN}$ with $S_0=0$, $S_1=\xi_1$, and  for $n\geq 2$, $S_n=S_{n-1}+\xi_n$
with $\xi_1$ uniformly distributed in $\BbS^2(0,1)$, and 
$\xi_n$ uniformly distributed in
 $\BbS^2(S_{n-1},1)\cap C(\xi_{n-1},S_{n-1})$, where $\BbS^2(x,r)$ is the
sphere of center $x$ and radius $r$.
Here the ``time'' $n$ indexing the process corresponds to the order
of appearance of a given monomer inside the macromolecular chain. For any
bounded measurable function $f:\BbR^3\rightarrow \BbR$, we have then
\[
\BbE(f(S_n)|\cF_{n-1})=\int_0^{2\pi} f(S_{n-1}+\sin\theta \cos\phi e_1+
\sin\theta\sin\phi e_2+\cos\theta e_3))\frac{d\phi}{2\pi},\]
where $\cF_n=\sigma(S_k, k\leq n)$ and $e_1,e_2,e_3$ is the canonical basis
of $\BbR^3$.
A natural simplification of the model consists in considering $(\xi_n)$ a
sequence of independent and identically distributed random variables,
getting thus a simple random walk on $\BbR^3$, an object that
has been extensively studied and generalised in various respects
and especially on non-commutative groups. It is not the intention of
the authors to report further in this direction since it is
perfectly well known by the community of probabilists and excellent
monographs have been devoted to the subject (see \cite{Spi}, \cite{Revuz},
\cite{Woe} for instance), but to report on some aspects 
developed mainly by physicists  and less known by probabilists.

The model of simple random walk is too na\"\i ve to realistically model
physical polymers: two different atoms cannot occupy the same
position. Hence a realistic model must be self-avoiding, spoiling
thus the Markovian character of the process.
Consider the simplest random walk on $\BbZ^d$, the nearest
neighbour random walk, \textit{i.e.}\ let
$E_d=\{\pm e_1, \ldots, \pm e_d\}$, where $(e_i)_{i=1,\ldots,d}$ denote
the standard basis of $\BbZ^d$, and let
$(\xi_n)_{n\in\BbN}$ be an independent and identically uniformly distributed
sequence of $E_d$-valued random variables. Then the process
defined by $S_0=0$ and $S_n=S_{n-1}+\xi_n$ for $n\geq 1$, provides
the Markovian description of the ordinary random walk. If we are
interested only to a finite sequence
$(S_n)_{n=0}^N$, an equivalent description is provided
by the trajectory space
\[\Omega_N=\{\omega:\{0,\ldots, N\}\rightarrow\BbZ^d\ \mid\
  \omega(0)=0, \omega(i)-\omega(i-1)\in E_d,
i=1,\ldots, N\},\]
equipped with the uniform probability measure
$\mu_N(\omega)=1/c_N$ for all $\omega\in \Omega_N$ where
$c_n=\card \Omega_N= (2d)^N$. 
For $k: 0\leq k\leq N$, the canonical projection
$S_k(\omega)=\omega_k$ has the same law as the random walk
defined by the sum $\sum_{i=1}^k\xi_i$, showing thus the equivalence
of the Markovian and trajectorial descriptions for simple random walks.
Adding the self-avoiding condition can be performed on the trajectorial
description but not on the Markovian one. More precisely, let
\[\Omega^\saw_N=\{\omega\in \Omega_N:\omega(i)\ne\omega(j), \ \textrm{for}\ 
0\leq i<j\leq N\}\]
and 
$\mu_N^\saw(\omega)=1/c_N^\saw$ for all $\omega\in\Omega_N^\saw$ with
$c_N^\saw=\card \Omega_N^\saw$. Notice however that the numerical
sequence $(c_N^\saw)_{N\in\BbN}$ is not explicitly known for $d\geq 2$ hence
the model of self-avoiding random walk has been so far intractable. Nevertheless, the
sequence of probability measures $(\mu^\saw_N)_{N\in\BbN}$ --- the probability
$\mu^\saw_N$ being defined on $\Omega_N^\saw$ for every $N\in\BbN$ --- is perfectly well
defined although intractable.

Instead of defining 
$\mu^\saw_N$ on $\Omega_N^\saw$
we can also define it on $\Omega_N$ by 
\begin{equation}
\label{eq:meas-saw}
\mu^\saw_N (d\omega)=\frac{1}{Z_N^\saw} \id_{\Omega_N^\saw}(\omega) 
\mu_N(d\omega),
\end{equation}
where $Z_N^\saw$ is a normalising factor, in fact $Z_N^\saw=c_N^\saw/(2d)^N$.
Physicists have introduced various approximations to deal with the untractable
measure $\mu^\saw_N$. One of them consists in approximating the indicator
appearing in the previous formula and defining
\begin{equation}
\label{eq:meas-edwards}
\mu_{N,\beta}(d\omega)=\frac{1}{Z_n(\beta)}\exp(-\beta H_N(\omega)) 
\mu_N(d\omega)
\end{equation}
where $H_N(\omega)=\card\{k: 2\leq k\leq N| \exists j<k \ \textrm{with}\
\omega(j)=\omega(k)\}$ (or some variant of this number)
counts the self-intersections of the trajectory $\omega\in
\Omega_N$, the real parameter $\beta$ is non-negative, and 
$Z_N(\beta)$ is a normalising factor, in fact $Z_N(\beta)=\int_{\Omega_N}
\exp(-\beta H_N(\omega))\mu_N(d\omega)$. The continuous version of this
model has been introduced by Edwards in \cite{Edw}. The discrete
version, defined in (\ref{eq:meas-edwards}), is known as discrete Edwards random walk or weakly self-avoiding
random walk. The reason for this nomenclature is the following.
For fixed $N$, we have that $\lim_{\beta\rightarrow 0}\mu_{N,\beta}=\mu_N$
and  $\lim_{\beta\rightarrow \infty}\mu_{N,\beta}=\mu^\saw_N$
so that the weakly self-avoiding walk interpolates between ordinary
and self-avoiding random walk. A much more difficult limit to study
is  $N\rightarrow \infty$ for fixed $\beta \in ]0,\infty[$; this limit
is dimension dependent and several authors have contributed
to its study 
\cite{Sym,BrySpe,Westwater-I,Westwater-I,Law-book,Law-anp,HarSla-saw5,MadSla,Guillotin-Zoladek,Vermet,Cadre,HofHolSla,KouPasPet,ForPasPet,Hue},
either rigorously or numerically\footnote{The day of submission of the present contribution,
we learnt about a new result on weakly self-avoiding walks \cite{BolRit}.}.

The formula \ref{eq:meas-edwards} has an interest far beyond its application
to the study of self-avoiding random walks since it is reminiscent of ideas
at the basis of Gibbs formulation of statistical mechanics and quantum
field theory. The Radon-Nikod\'ym derivative 
$d\mu_{N,\beta}/d\mu_N=\exp(-\beta H_N(\omega))/Z_N(\beta)$ is interpreted
as a Boltzmann factor, making more rare  trajectories with many self-intersections
in the statistical sample described by the measure $\mu_N$. The
limit $N\rightarrow \infty$ is also a standard procedure in statistical
mechanics and quantum field theory known as thermodynamic limit
or infrared limit respectively. Thus, this formula
is naturally connecting random walks and various other physical theories
defined on an infinite graph for which the fundamental quantities
can be written as (formal) random walk expansions. In some particular
situations, these formal expansions converge; they provide therefore
an valuable probabilistic tool for the study of the asymptotic
behaviour of physical models.

\subsection{An example of random walk expansion}
It is immediate to see that the Markov operator for a simple
random walk on an undirected graph is essentially the discrete
Laplacian on the graph.
Random walk expansions can take a very sophisticated formulation; all
of them can be seen however as (non-trivial) generalisations of a very
simple formula of inversion of the Markov operator that 
connects the Green's function of the regularised Laplacian with
a power series on random walks of arbitrary length given in the following

\begin{lemma}
Let $\Delta$ be the difference Laplacian on $\BbZ^d$ and $m\ne 0$ a fixed
parameter (the free mass of a quantum field theory). Then for any $x,y\in\BbZ^d$, 
\[(m^2\id -\Delta)^{-1}_{xy}=\sum_{\omega\in G^*(x,y)} (2d+m^2)^{-|\omega|},\]
where $G^*(x,y)=\cup_{n=0}^\infty G^n(x,y)$ and
the set $G^n(x,y)$ is the set of paths of length $n$ on the non-oriented
graph $\BbG$ having  $\BbZ^d$ as vertex set and
with nearest neighbours vertices as edge set,  that start at point $x$ and end
at point $y$.
\end{lemma}

\proof 
Write $\Delta=J-2d\id$ where 
\[J_{xy}=\left\{\begin{array}{ll}
1 & \textrm{if} \  y-x\in E_d\\
0 & \textrm{otherwise}
\end{array}\right.\]
and develop $((m^2+2d)\id-J)^{-1}$ as formal series in powers of $J$.
For $m\ne 0$ the series converges defining thus the left-hand side of the formula. 
\qed

Consider now a general graph $\BbG=(G^0,G^1,r,s)$ with $G^1$ a particular
subset of $G^0\times G^0\setminus \diag (G^0\times G^0)$ (\textit{i.e.}\
the graph is simple and without loops); we assume moreover that the graph
is undirected, \textit{i.e.} if $(u,v)\in  G^1$ then $(v,u)\in G^1$. 
Let $J$ and $L$ be $G^0\times G^0$ matrices
such that $J_{uv}>0$ if $v\in N_u$ and
$J_{uv}=0$ otherwise and $L_{uv}= \lambda_u \delta_{uv}$,
with $\lambda_u>0$ for all $u\in G^0$.
\begin{lemma}
Suppose that the parameters $(\lambda_u)_{u\in G^0}$ are large enough. Then
\begin{enumerate}
\item
\[(L- J)^{-1}_{uv}= \sum_{\omega\in G^*(u,v)}
\prod_{a\in \omega} J_a \prod_{v\in G^0} \lambda_v^{-\eta_{|\omega|}(v,\omega)},\]
where $\eta_{|\omega|}(v,\omega)= \sum_{k=0}^{|\omega|} \id_{\{v\}} (\omega_k)$ is the
occupation time of the vertex $v$ by the trajectory $\omega$.
\item
\[\det (L- J)^{-1} = \left(\prod_{v\in G^0} \lambda_v\right)^{-1}
\exp\left(\sum_{\omega\in \cL^*} J_\omega 
(\prod_{u\in G^0} \lambda_u^{-\eta_{|\omega|}(u,\omega)}\right),\]
where $J_\omega=\prod_{a\in\omega} J_a$ and $\cL^*$ is the set of loops of arbitrary length
(\textit{i.e.}\ equivalence classes of random walks of arbitrary length starting
and ending on the same vertex.) 
\end{enumerate}
\end{lemma}
\proof
The matrix $L$ is diagonal and invertible. The formulae are
again obtained by standard formal power series expansions that converge
when $\lambda_u$ are chosen sufficiently large (see \cite{BryFroSpe} for details.) 
\qed

Consider always a simple undirected graph without loops
$\BbG=(G^0,G^1,r,s)$ and let $f_v:\BbR\rightarrow\BbR$ be a family
of maps indexed by $v\in G^0$ such that 
$\lim_{t\rightarrow\infty}f_v(t)\exp(ct) =0$ for some $c>0$.
Denote $\BbX=\{x:G^0\rightarrow \BbR^\nu\}\simeq (\BbR^\nu)^{G^0}$ the
configuration space, so that each configuration $x\in\BbX$ is
the collection $(x_v)_{v\in G^0}$ with $x_v\in\BbR^\nu$. The space
$ \BbR^\nu$ is equipped with its Borel $\sigma$-algebra,
$\cB(\BbR^\nu)$ and let $\kappa_u$
be a family of  continuous measures on $(\BbR^\nu, \cB(\BbR^\nu))$ 
defined to have
a density with respect to the $\nu$-dimensional Lebesgue measure
$\kappa_u(dx_u)=f_u(x^2)dx_u=f_u(x_{u,1}^2+\ldots x_{u,\nu}^2)
dx_{u,1}\ldots dx_{u,\nu}$. The configuration space can be naturally
equipped with a product measure structure $\prod_{u\in G^0} \kappa_u$.
However, a product measure structure is not very interesting for physical
purposes since it corresponds to an infinite system of non-interacting
components. To introduce some interaction, let $J$ be an infinite
$G^0\times G^0$ matrix with $J_{uv}=J_{vu}>0$ when $v\in N_u$ and 
$J_{uv}=0$ otherwise. Since the matrix is symmetric, it defines a quadratic
form on $\BbX$, known as (formal) Hamiltonian
$H(x)=-\frac{1}{2}\sum_{u,v\in G^0}( x_u, J_{uv} x_v)\equiv
-\frac{1}{2}\sum_{u,v\in G^0} \sum_{\alpha=1}^\nu x_{u,\alpha} J_{uv} 
x_{v,\alpha}$. We can then define (formally) a probability $\mu$
on $(\BbX,\cF)$, where $\cF$ is the natural $\sigma$-algebra, 
by 
\[\mu(A)=\frac{1}{Z} \int_A \exp(-H(x)) \prod_{v\in G^0} (f_v(x_v^2) dx_v)),\]
where $Z$ is a normalising factor.
There are various standard procedures to give a mathematical meaning to the
above expressions. One of them is the following: suppose that the
graph $\BbG$ is isometrically embedded in $\BbR^d$ for some $d$. Let
$\Lambda_n=[-n,n]^d\cap G^0$ be the set of vertices of the graph
inside an hypercubic box of size  $2n+1$. Then, we define
the finite volume Hamiltonian
\[H_n(x)=-\frac{1}{2} \sum_{u,v\in\Lambda_n} (x_u, J_{uv} x_v)\]
and the finite volume probability
\[\mu_n(A)= \frac{1}{Z_n}\int_A  
\exp(-H_n(x)) \prod_{v\in \Lambda_n} (f_v(x_v^2) dx_v)),\]
 where $Z_n$ is the corresponding normalising factor.
The sequence $(\mu_n)$ is a perfectly well defined sequence of probability
measures. When this sequence converges weakly, we call the weak limit
infinite volume Gibbs measure associated with the Hamiltonian $H$.
Notice however that the existence of the weak limit is highly non
trivial and it is granted only for some cases 
(see  \cite{Sym,BryFroSpe,GliJaf,Mal} for instance). 

\begin{theorem}[Symanzik]
Assume that the Hamiltonian is such that the weak limit $\mu$ exists. 
For $p=1,2,\ldots$ let $\{v_1,\ldots , v_{2p}\}$ be a given set of vertices;
partition this set into $p$ disjoint pairs. For each such pair of vertices
let
$\omega^{(l)}$ be a random walk of arbitrary length  starting in one vertex
and ending to the other vertex of the pair.
Then 
\begin{eqnarray*}
\mu(x_{v_1,\alpha_1}\cdots x_{v_{2p}, \alpha_{2p}})&\equiv&
\int x_{v_1,\alpha_1}\cdots x_{v_{2p}, \alpha_{2p}} \mu(dx)\\
&=& \sum_{\omega^{(1)}, \ldots, \omega^{(p)}} \frac{W(\omega^{(1)},
\ldots, \omega^{(p)})}{Z},
\end{eqnarray*}
where the sum extends over all partitions of the
set of vertices and all random walk defined on the $p$ pairs
and
\begin{eqnarray*}
W(\omega^{(1)},\ldots, \omega^{(p)})&=&\sum_{n=0}^\infty
\frac{1}{n!} (\frac{\nu}{2})^n \sum_{\gamma ^{(1)},\ldots,\gamma^{(n)}\in \cL^*}
\prod_{m=1}^n J_{\gamma^{(m)}} \prod_{l=1}^p J_{\omega^{(l)}}\\
&&
\exp(-U(\gamma ^{(1)},\ldots,\gamma^{(n)}, \omega^{(1)},\ldots, \omega^{(p)})),\end{eqnarray*}
where $\cL^*$ denotes the set of loops of arbitrary length.
Moreover, 
the normalising factor reads
\[Z= \sum_{n=0}^\infty
\frac{1}{n!} (\frac{\nu}{2})^n \sum_{\gamma^{(1)},\ldots,\gamma^{(n)}\in \cL^*}
\prod_{m=1}^n J_{\gamma^{(m)}} 
\exp(-U(\gamma^{(1)},\ldots,\gamma^{(n)}))\]
and the mapping $U$ defined on an arbitrary Cartesian product of random
walks is given by
\[\exp(-U(\omega^{(1)},\ldots, \omega^{(k)}))=
\prod_v\int_\Gamma \hat{f}_v(z_v) 
(2\i z_v)^{-h_v(\omega^{(1)},\ldots, \omega^{(k)})} dz_v,\]
$\hat{f}$ being the Fourier transform of $f_v$,
 $h_v(\omega^{(1)},\ldots, \omega^{(k)})= \sum_{l=1}^k
\eta_{|\omega^{(l)}|}(v,\omega^{(l)})+\frac{\nu}{2}$, and $\Gamma$ an appropriately
chosen integration contour.
\end{theorem}

\Rk The previous theorem looks formidable. It is worth noticing however
that it is nothing more than a clever combinatorial recombination
of terms appearing in the power series expansion of the exponential and that
far-reaching results in quantum field theory \cite{AizFro,FerFroSok}
--- impossible or much more difficult 
to obtain otherwise --- are obtained
by the random walk representation it provides. Variants of this representation --- known
under the generic name of cluster expansions or abstract polymer models \cite{Mal,GliJaf} ---
are used in many different contexts, like statistical mechanics, disordered
systems etc.\ and in the cases it can be rigorously applied it provides
a very powerful probabilistic tool for the study of covariance properties
of limit Gibbs measures. 

\section{The physical relevance of random walks on oriented lattices}
\label{sec-rwol}
For any simple directed graph without loops $\BbG=(G^0,G^1,r,s)$ 
we define two operators \cite{Paschke},
$D:l^2(G^0)\rightarrow l^2(G^1)$ by $(Df)(a)=f(s(a))-f(r(a))$ for every 
$a\in G^1$ and its adjoint 
 $D^*:l^2(G^1)\rightarrow l^2(G^0)$ by 
$(D^*\phi)(v)=\sum_{a\in s^{-1}(v)}\phi(a)-
\sum_{a\in r^{-1}(v)}\phi(a)$ for every 
$v\in G^0$. Then $(-D^*Df)(v)= -d_v f(v)+ \sum_{u\in N_v} f(u)\equiv
(\Delta f)(v)$.
For a simple random walk on $\BbG$, defined by its stochastic matrix 
$(P(u,v))_{u,v\in G^0}$, the Markov operator $M$ defined on bounded
functions $f$ by $Mf(u)=\sum_{u\in G^0} P(u,v)f(v)- f(u)$ --- contrary to the
case of unoriented graphs where it is expressed in terms of the Laplacian ---
cannot be expressed in terms of the Laplacian since
$Mf(u)=\frac{1}{d_u^-}\sum_{a\in s^{-1}(v)} Df(a)
\ne \frac{1}{d_u^-}(\Delta f)(u)$.
As a matter of fact, $M$ is (roughly) reminiscent of the Dirac operator
on the graph, providing thus the first hint that random walks on oriented
lattices are relevant for non-commutative geometry.

\subsection{A $C^*$-algebraic description of oriented lattices}
\label{ssec-Cstar}
With every oriented graph we can associate a $C^*$-algebra of operators, known
as the Cuntz-Krieger algebra \cite{CunKri} of the graph 
\cite{KumPasRae,KumPasRen}.

Let $(V_i)_{i\in I}$ be a finite or denumerable family of non-zero partial
isometries and $A$ a $I\times I$ $\{0,1\}$-valued matrix whose rows contain
a finite number of ones. The Cuntz-Krieger algebra, $\cO_A$, associated
with the matrix $A$ is the $C^*$-algebra, defined up to isomorphisms, by
the relations
\[V^*_i V_i=\sum_{j\in I} A_{ij} V_j V^*_j, \  i\in I.\]
The connection of the $\cO_A$ algebra with oriented graphs is done as follows.
Let $\BbG=(G^0,G^1,r,s)$ be a row-finite, locally finite graph and consider
the corresponding path space $G^*$ defined in section \ref{sec-notations}.
Let $(P_v)_{v\in G^0}$ be a set of mutually orthogonal projections
and $(V_a)_{a\in G^1}$ a set of non-zero partial isometries satisfying
\[V^*_a V_a =P_{r(a)}, \forall a\in G^1 \  \textrm{and}\
P_v=\sum_{a\in s^{-1}(v)} V_a V^*_a, \forall v\in G^0.\]
Define the edge-matrix $(A_\BbG(a,b))_{a,b \in G^1}$ of the graph $\BbG$
by
\[A_\BbG(a,b)=\left\{\begin{array}{ll}
1 & \textrm{if} \ r(a)=s(b)\\
0 & \textrm{otherwise.}
\end{array}\right.\]
Then 
\[P_{r(a)}=V^*_a V_a=\sum_{b\in G^1: r(a)=s(b)} V_b V^*_b=
\sum_{b\in G^1} A_\BbG(a,b) V_b V^*_b.\]
The $G^1$ indexed Cuntz-Krieger $C^*$-algebra $\cO_{A_\BbG}$ is called
the $C^*$-algebra of the graph $\BbG$ and is denoted $C^*(\BbG)$. The partial
isometries $V_a$ defined on $G^1$ are naturally extended on the path space
$G^*$ for every $\alpha\in G^*$ by
\[ V_\alpha= V_{a_1}\cdots V_{a_{|\alpha|}}.\]
Then we have \cite{KumPasRae} the
\begin{theorem}[Kumjian, Pask, and Raeburn]
Let $\{P_v, V_a; v\in G^0, a\in G^1\}$ be the Cuntz-Krieger algebra
associated 
with $\BbG$ and let $\beta$ and $\gamma$ be arbitrary paths of $G^*$.
Then
\[V^*_\beta V_\gamma=\left\{\begin{array}{ll}
V_{\gamma'} & \textrm{if}\ \gamma= \beta \gamma', \gamma'\ni G^0\\
P_{r(\gamma)} & \textrm{if}\ \gamma= \beta\\
V^*_{\beta'} & \textrm{if}\ \beta=\gamma \beta', \beta'\ne G^0\\
0 & \textrm{otherwise.}
\end{array}\right.\]
\end{theorem}
This construction associates with every path in $G^*$ an operator
of $C^*(\BbG)$. Therefore the $C^*$-algebra $C^*(\BbG)$ can be thought as
the non-commutative analogue of the subshift space 
\cite{BraJorPal-iterated}, corresponding
to the matrix $A_\BbG$.
It is worth noticing that beyond directed lattices, the Cuntz-Krieger algebras
are ultimately connected with various other topics, like wavelets \cite{BraJorPal-representation},
tilings \cite{Kellendonk,BelBenGam}, generalised sub-shifts \cite{BraJorPal-iterated}, non-commutative
geometry \cite{Connes-book}, etc. 

\subsection{Non-reversible evolution of quantum states}
\label{ssec-evolution}
It is convincingly argued lately that the progress in semiconductor
technology will reach the limits of applicability of classical reasoning
used in information theory and computer science because quantum effects
will start to be determining \cite{WilCle,Preskill,Holevo,Pet-QI}.

In classical physics a microscopic state of a multi-component system 
is described as an element of the Cartesian product of state spaces
of individual components; for instance to determine
the microstate of a litre of gas, it is necessary to know the precise positions
and momenta of all its molecules. The set of all 
microstates is called \textit{configuration space}. Macroscopic states
are classically probability measures on the configuration space and
physical observables are bounded measurable real-valued functions
on the configuration space. Time evolution is implemented as
a Markov semi-group acting on macrostates; it can be reversible
when the stochastic kernel of the semigroup is deterministic, or
irreversible in general.

In quantum physics the configuration space of microstates is a complex
separable Hilbert space $\cH$ (in general infinite dimensional but finite
dimensional spaces cases or also of interest.) The macrostates
are self-adjoint, positive operators with  unit trace  that are projections.
Such operators are usually called \textit{density matrices}. Observables
are self-adjoint bounded operators acting on $\cH$. 
Time evolution is implemented
by a completely positive  transformation on the set of density matrices 
$\rho\rightarrow \sum_{i\in I}T_i \rho T_i^*$. Such evolutions can
be reversible --- when $I=\{0\}$ is a singleton and $T_0$ is unitary ---
or irreversible --- when we only require that $\sum_{i\in I} T_i T^*_i\leq 
\mathrm{id}_\cH$.  One immediately remarks that $C^*$-algebras provide a unified
approach to both classical (Abelian ones) and quantum (non Abelian ones) systems.

In the context of applications in quantum 
information and communication, unitary
transformations correspond to quantum logical gates while irreversible ones
to noisy transmission through quantum channels or to measurements. 
When a quantum macrostate
$\rho$ is transmitted through a noisy channel, different operators
from the family $(T_i)_{i\in I}$ will sequentially act on $\rho$ to get
$\rho\rightarrow \hat{\rho}_n =\sum_{i_1,\ldots,i_n} T_{i_n}\cdots T_{i_1} \rho
T_{i_1}^* \cdots T_{i_n}^*$.
It is thus clear that products of operators will appear in the form
of  products of partial isometries
along paths of oriented graphs. The main difference is that
we have not required here the evolution to be implemented
by partial isometries but by general non-commuting operators
satisfying only $\sum_{i\in I} T_i T^*_i\leq \mathrm{id}_\cH$.
Nevertheless the analogy can me made complete by virtue of the following
result \cite{Pop}
\begin{theorem}[Popescu]
For every sequence (finite or denumerable) $(T_i)_{i\in I}$ of
non-commuting operators acting on a Hilbert space $\cH$, such
that $\sum_{i\in I} T_i T^*_i\leq \mathrm{id}_\cH$, there exists a minimal
isometric dilation into operators $(V_i)_{i\in I}$  acting on a Hilbert
space $\cK\supset \cH$, uniquely determined up to isomorphisms, such that
the family $(V_i)_{i\in I}$ is composed by non-zero partial isometries
on $\cK$.
\end{theorem}

Now it is evident that random walks on a (partially) directed lattice
induce a random walk on the space of density matrices; recovering
information from the perturbed macrostate $\hat{\rho}_n$ will be
better when $\hat{\rho}_n$ is close (with respect to an appropriate topology)
to $\rho$. This remark is made for the first time in \cite{Ler} and gives
a practical physical and technological relevance to questions
of recurrence of the random walk on a randomly oriented lattices.

\section{Sketch of the proofs}
\label{sec-proofs}
The idea of the proof of theorems \ref{th-L}, \ref{th-H}, and
\ref{th-O} is to decompose the Markov chain $(\bM_n)$ into a
vertical skeleton $(Y_n)$ and an horizontal embedded random walk $(X_n)$ that
--- when sampled on a particular sequence of random times defined in terms
of the vertical skeleton --- has the same recurrence/transience properties
as the original random walk $(\bM_n)$.

Let $(\psi_n)_{n\in\BbN_+}$ be a sequence of independent, identically
distributed, $\{-1,1\}$-valued symmetric Bernoulli variables and define $Y_0=0$ and
\[Y_n=\sum_{k=1}^n\psi_k, \ \ n=1,2,\ldots,\]
the simple $\BbV_2$-valued symmetric one-dimensional random walk.
We call the process $(Y_n)_{n\in\BbN}$ the \textit{vertical skeleton}.
We denote by
\[\eta_n(y)=\sum_{k=0}^n\id_{\{Y_k=y\}}, \ \ n\in\BbN, y\in\BbV_2\]
its \textit{occupation time} at level $y$.
If $\cF_n=\sigma(\psi_i, i\leq n)$ then $\eta_n(y)$ is obviously 
$\cF_n$-measurable.

Define $\sigma_0=0$ and recursively, for $n=1,2,\ldots$,
$\sigma_n=\inf\{k\geq\sigma_{n-1}: Y_k=0\}>\sigma_{n-1}$,
the $n^{\textrm{th}}$ \textit{return to the origin} for the vertical
skeleton. It is known by the standard theory of simple symmetric
one-dimensional random walk that almost surely
$\sigma_n<\infty$ for all $n$.

To define the horizontal embedded random walk, let
$(\xi^{(y)}_n)_{n\in\BbN_+, y\in\BbV_2}$ be a doubly infinite sequence
of independent identically distributed $\BbN$-valued
geometric random variables of parameters $p=2/3$ and $q=1-p=1/3$, 
\textit{i.e.}\
$\BbP(\xi^y_1=\ell)=pq^\ell, \ \ \ell=0,1,2,\ldots$.
Let  moreover
\[T_n=n+\sum_{y\in\BbV_2} \sum_{i=1}^{\eta_{n-1}(y)} \xi_i^{(y)}\]
be the instant just after the random walk $(\bM_k)$ has performed
its $n^{\textrm{th}}$ vertical move (with the convention that the
sum $\sum_{i=1}^{\eta_{n-1}(y)}$ vanishes whenever $\eta_{n-1}(y)=0$.) Then
\[\bM_{T_n}=(X_n,Y_n),\]
where $(Y_n)$ is the vertical skeleton and 
\[X_n= \sum_{y\in\BbV_2}\epsilon_y \sum_{i=1}^{\eta_{n-1}(y)} \xi_i^{(y)}\]
represents the total horizontal displacement when the random walk
$(\bM_k)$ has completed exactly $n$ vertical moves. Notice also
that the horizontal embedded random walk $(X_n)$ can be viewed as a
random walk with unbounded jumps in a random scenery and
that $\bM_{T_{\sigma_n}}=(X_{\sigma_n},0)$. Now
$(\bM_n)$ can return to the origin if and only if both vertical and horizontal
components vanish simultaneously. Since the vertical component can vanish
only at instants $\sigma_n, n\in\BbN$, we prove in \cite{CamPet-rwrol}
the following 
\begin{lemma}
\begin{enumerate}
\item  $\sum_{n=0}^\infty
\BbP(X_{\sigma_n}=0)=\infty$ if and only if the random walk 
$(\bM_l)$ is recurrent.
\item 
$\sum_{n=0}^\infty
\BbP(X_{\sigma_n}=0)<\infty$ if and only if  the random walk
$(\bM_l)$ is transient.
\end{enumerate}
\end{lemma}

Introduce now the characteristic function $\chi(\theta)=
\BbE\exp(i\theta\xi_1^{(0)})= \frac{p}{1-q\exp(i\theta)}=
r(\theta)\exp(i\alpha(\theta))$ and observe that $r$ is an even
function of $\theta$ while $\alpha$ is odd. Hence, denoting
$\cF=\cF_\infty$ and $\cG=\sigma(\epsilon_y, y\in\BbV_2)$, we have
\begin{eqnarray*}
\BbE\exp(i\theta X_n) 
&=&\BbE\left(\BbE(\exp(i\theta\sum_{y\in\BbV_2}\epsilon_y
\sum_{i=1}^{\eta_{n-1}(y)}\xi_i^{(y)}|\cF\vee\cG)\right)\\
&=&\BbE\left(\prod_{y\in\BbV_2}\chi(\theta\epsilon_y)^{\eta_{n-1}(y)}\right)\\
&=&\BbE(r(\theta)^{\sum_{y\in\BbV_2}\eta_{n-1}(y)}
\exp(i\alpha(\theta)\Delta_n))\\
&=& r(\theta)^n \BbE(\exp(i\alpha(\theta)\Delta_n)),
\end{eqnarray*}
where $\Delta_n=\sum_{y\in\BbV_2}\epsilon_y\eta_{n-1}(y)$.

\begin{lemma}
For the $\BbL$ lattice, $\Delta_{\sigma_n}=0$.
\end{lemma}

\proof This is an elementary combinatorial result whose complete
proof is given in \cite{CamPet-rwrol}.
\qed

The proof of theorem \ref{th-L} follows now immediately since
$\BbE\exp(i\theta X_{\sigma_n})=\BbE(r(\theta)^{\sigma_n})=
(\sqrt{1-r(\theta)^2})^n$. Hence
$\BbP(X_{\sigma_n}=0)= \lim\limits_{\epsilon\rightarrow 0}
2 \int_\epsilon^{\pi} \frac{1}{\sqrt{1-r(\theta)^2}} d\theta=\infty$.

\Rk Notice that for the random walk on the $\BbL$ lattice, various other more
elegant proofs can be given; we presented the most elementary one.

For the corresponding result on the $\BbH$ lattice, the proof is based
on the following lemma shown in \cite{CamPet-rwrol}
\begin{lemma}
Denote by $(\rho_k)_{k\in\BbN}$ a sequence of independent
identically distributed Rademacher variables and $(\tau_k)_{k\in\BbN}$
a sequence of independent, identically distributed random variables,
independent of the sequence $(\rho_k)_{k\in\BbN}$, such that
$\tau_1\elaw \sigma_1$, \textit{i.e.}\ the random variables $\tau_k$
have the same law as the time of first return to the origin for the
skeleton random walk. Then
\[\Delta_{\sigma_n} \elaw
\sum_{k=1}^n \rho_k (\tau_k-1) +n.\]
\end{lemma}

The proof of the theorem \ref{th-H} follows then from the equality
$  \BbE\left(\exp(i\theta X_{\sigma_n})\right)=
g(\theta)^n$, 
where $g(\theta)=\frac{1}{2}\chi(\theta)
\left[ \left(1-\sqrt{1-\chi(\theta)^2}\right)\exp(-i\alpha(\theta)) +
\left(1-\sqrt{1-\ol{\chi}(\theta)^2}\right)\exp(i\alpha(\theta)) \right]$.
This expression allows an explicit estimate for $\sum_{n=1}^\infty
\BbP(X_{\sigma_n}=0)$ that converges in the present case, 
establishing transience of the random 
walk.

Finally for the $\BbO_{\bom{\epsilon}}$ lattice, no simple decomposition
can be made and we need joint estimates. The idea of the proof is to decompose
the probability $p_n=\BbP(X_{2n}=0, Y_{2n}=0)$ into
$p_n=p_{n,1}+p_{n,2}+p_{n,3}$,
where  
\begin{eqnarray*}
p_{n,1}&=& \BbP(I(X_{2n},-\epsilon_0 Z)\ni 0; Y_{2n}=0;B_n)\\
p_{n,2}&=& \BbP(I(X_{2n},-\epsilon_0 Z)\ni 0; Y_{2n}=0;A_n\setminus B_n)\\
p_{n,3}&=& \BbP(I(X_{2n},-\epsilon_0 Z)\ni 0; Y_{2n}=0;A_n^c),
\end{eqnarray*}  
and $A_n=A_{n,1}\cap A_{n,2}$ with
\begin{eqnarray*}
A_{n,1} &=& \{\omega\in\Omega:
\max_{0\leq k\leq 2n}|Y_k|<n^{\frac{1}{2}+\delta_1}\}\
\ \textrm{for some}\ \ \delta_1>0,\\
A_{n,2} &=& \{\omega\in\Omega: \max_{y\in\BbV_2}\eta_{2n-1}(y)<
n^{\frac{1}{2}+\delta_2}\}\ \ \textrm{for some}\ \ \delta_2>0,\\
B_n &=& \{\omega\in A_n:
 \left|\sum_{y\in\BbV_2}\epsilon_y \eta_{2n-1}(y)\right|
>n^{\frac{1}{2}+\delta_3}\}\ \ \textrm{for some}\ \ \delta_3>0
\end{eqnarray*}
for some $\delta_1>0, \delta_2>0$ and $\delta_3>0$.
The technical part of the proof consists in showing that 
$p_{n,1}$ and $p_{n,3}$ are of order $\cO(\exp(-n^\delta))$ with 
some $\delta>0$ for large $n$ while the main part of the mass charging
the event $\{X_{2n}=0, Y_{2n}=0\}$ is supported by the set
$A_n\setminus B_n$. More precisely it is shown in \cite{CamPet-rwrol} that
on the set $A_n\setminus B_n$ we have
$\BbP(X_{2n}=0|\cF\vee\cG)=\cO(\sqrt{\frac{\ln n}{n}})$
while $\BbP(A_n\setminus B_n|\cF\vee\cG)=\cO(n^{-1/4+\delta_4})$
so that together with the standard estimate
$\BbP(Y_{2n}=0)=\cO(1/\sqrt{n})$ the series $\sum_n p_n$ is shown to converge.

Notice also that some additional results concerning the mean speed
and law of large numbers are presented in \cite{CamPet-rwrol}.

\bibliographystyle{plain}
\scriptsize{
\bibliography{rwre,saw,petritis,cstar,quinfo}

\begin{thebibliography}{10}

\bibitem{AizFro}
Michael Aizenman and J{\"u}rg Fr{\"o}hlich.
\newblock Topological anomalies in the $n$ dependence of the $n$-states {P}otts
  lattice gauge theory.
\newblock {\em Nuclear Phys. B}, 235(1, FS11):1--18, 1984.

\bibitem{BelBenGam}
Jean Bellissard, Riccardo Benedetti, and Jean-Marc Gambaudo.
\newblock {Spaces of tilings, finite telescopic approximations and
  gap-labelling}, preprint 2001, eprint {arXiv:math.DS/0109062}.

\bibitem{BolRit}
Erwin Bolthausen and Christine Ritzmann.
\newblock A central limit theorem for convolution equations and weakly
  self-avoiding walks.

\bibitem{BraJorPal-representation}
Ola Bratteli and Palle E.~T. Jorgensen.
\newblock A connection between multiresolution wavelet theory of scale ${N}$
  and representations of the {C}untz algebra ${\cO}\sb {N}$.
\newblock In {\em Operator algebras and quantum field theory (Rome, 1996)},
  pages 151--163. Internat. Press, Cambridge, MA, 1997.

\bibitem{BraJorPal-iterated}
Ola Bratteli and Palle E.~T. Jorgensen.
\newblock Iterated function systems and permutation representations of the
  {C}untz algebra.
\newblock {\em Mem. Amer. Math. Soc.}, 139(663):x+89, 1999.

\bibitem{BryFroSpe}
D.~Brydges, J.~Fr{\"o}hlich, and Th. Spencer.
\newblock The random walk representation of classical spin systems and
  correlation inequalities.
\newblock {\em Comm. Math. Phys.}, 83(1):123--150, 1982.

\bibitem{BrySpe}
David Brydges and Thomas Spencer.
\newblock Self-avoiding walk in $5$ or more dimensions.
\newblock {\em Comm. Math. Phys.}, 97(1-2):125--148, 1985.

\bibitem{Cadre}
B.~Cadre.
\newblock Une preuve standard du principe d'invariance de {S}toll.
\newblock In {\em S\'eminaire de Probabilit\'es, XXXI}, pages 85--102.
  Springer, Berlin, 1997.

\bibitem{CamPet-rwrol}
M.~Campanino and D.~Petritis.
\newblock {Random walks on randomly oriented lattices}, preprint 2001 submitted
  for publication, eprint arXiv:math.PR/0111305.

\bibitem{Connes-book}
Alain Connes.
\newblock {\em Noncommutative geometry}.
\newblock Academic Press Inc., San Diego, CA, 1994.

\bibitem{CunKri}
J.~Cuntz and W.~Krieger.
\newblock A class of ${C}\sp{\ast} $-algebras and topological {M}arkov chains.
\newblock {\em Invent. Math.}, 56(3):251--268, 1980.

\bibitem{ForPasPet}
Ph. de~Forcrand, J.~Pasche, and D.~Petritis.
\newblock Critical behaviour of {E}dwards random walk in two dimensions: a case
  where the fractal and {H}ausdorff dimensions are not equal.
\newblock {\em J. Phys. A}, 21(19):3771--3782, 1988.

\bibitem{Edw}
S.~F. Edwards.
\newblock The statistical mechanics of polymers with excluded volume.
\newblock {\em Proc. Phys. Soc.}, 85:613--624, 1965.

\bibitem{Ein}
A.~Einstein.
\newblock {Ü}ber die {T}heorie der {B}rownschewegung.
\newblock {\em Ann. der Physik}, 19:371, 1906.

\bibitem{FerFroSok}
R.\ Fern\'andez, J.\ Fr\"ohlich, and A.\ Sokal.
\newblock {\em Random walks, critical phenomena and triviality in quantum field
  theory}.
\newblock Springer-Verlag, Berlin, 1992.

\bibitem{Flo}
P.~Flory.
\newblock The configuration of a real polymer chain, 1949.

\bibitem{GliJaf}
J.~Glimm and A.~Jaffe.
\newblock {\em Quantum physics}.
\newblock Springer-Verlag, New York, second edition, 1987.
\newblock A functional integral point of view.

\bibitem{Guillotin-Zoladek}
N.~Guillotin.
\newblock Marches al\'eatoires \`a une dimension avec auto-int\'eraction,
  d'apr\`es {H}. {Z}oladek.
\newblock In {\em S\'eminaires de Probabilit\'es de Rennes (1995)}, page~12.
  Univ. Rennes I, Rennes, 1995.

\bibitem{HarSla-saw5}
Takashi Hara and Gordon Slade.
\newblock Self-avoiding walk in five or more dimensions. {I}. {T}he critical
  behaviour.
\newblock {\em Comm. Math. Phys.}, 147(1):101--136, 1992.

\bibitem{Holevo}
A.~Holevo.
\newblock {\em Statistical structure of quantum theory}.
\newblock Springer-Verlag, Berlin, 2001.

\bibitem{Hue}
I.~Hueter.
\newblock Proof of the conjecture that planar self-avoiding walk has root mean
  square displacement exponent 3/4, 2001.

\bibitem{Kellendonk}
Johannes Kellendonk.
\newblock Noncommutative geometry of tilings and gap labelling.
\newblock {\em Rev. Math. Phys.}, 7(7):1133--1180, 1995.

\bibitem{KouPasPet}
F.~Koukiou, J.~Pasche, and D.~Petritis.
\newblock The {H}ausdorff dimension of the two-dimensional {E}dwards' random
  walk.
\newblock {\em J. Phys. A}, 22(9):1385--1391, 1989.

\bibitem{KumPasRae}
A.~Kumjian, D.~Pask, and I.~Raeburn.
\newblock Cuntz-{K}rieger algebras of directed graphs.
\newblock {\em Pacific J. Math.}, 184(1):161--174, 1998.

\bibitem{KumPasRen}
A.~Kumjian, D.~Pask, I.~Raeburn, and J.~Renault.
\newblock Graphs, groupoids, and {C}untz-{K}rieger algebras.
\newblock {\em J. Funct. Anal.}, 144(2):505--541, 1997.

\bibitem{Law-anp}
G.~F. Lawler.
\newblock The infinite self-avoiding walk in high dimensions.
\newblock {\em Ann. Probab.}, 17(4):1367--1376, 1989.

\bibitem{Law-book}
G.~F. Lawler.
\newblock {\em Intersections of random walks}.
\newblock Birkh\"auser Boston Inc., Boston, MA, 1991.

\bibitem{Ler}
Ph. Leroux.
\newblock Phd thesis, {U}niversit\'e de {R}ennes {I}, in preparation 2002.

\bibitem{MadSla}
N.~Madras and G.~Slade.
\newblock {\em The self-avoiding walk}.
\newblock Birkh\"auser Boston Inc., Boston, MA, 1993.

\bibitem{Mal}
V.~A. Maly{\v{s}}ev.
\newblock Cluster expansions in lattice models of statistical physics and
  quantum field theory.
\newblock {\em Uspekhi Mat. Nauk}, 35(2(212)):3--53, 279, 1980.

\bibitem{MenPet-wedge}
M.\ Menshikov and D.\ Petritis.
\newblock Markov chains in a wedge with excitable boundaries, preprint
  université de {R}ennes 1, 2000.

\bibitem{Pai}
A.\ Pais.
\newblock {\em Subtle is the Lord\ldots}.
\newblock Oxford University Press, Oxford, 1982.

\bibitem{Paschke}
W.~L. Paschke.
\newblock The flow space of a directed ${G}$-graph.
\newblock {\em Pacific J. Math.}, 159(1):127--138, 1993.

\bibitem{Pet-QI}
D.~Petritis.
\newblock Introduction to quantum information and communication, in preparation
  2002.

\bibitem{Pol}
G.\ P\'olya.
\newblock {Ü}ber eine {A}ufgabe der {W}ahrscheinlichkeitsrechnung betreffend
  die {I}rrfahrt in {S}tra{\ss}ennetz.
\newblock {\em Math. Ann.}, 84:149--160, 1921.

\bibitem{Pop}
G.~Popescu.
\newblock Isometric dilations for infinite sequences of noncommuting operators.
\newblock {\em Trans. Amer. Math. Soc.}, 316(2):523--536, 1989.

\bibitem{Preskill}
J.~Preskill.
\newblock Quantum computing: pro and con.
\newblock {\em R. Soc. Lond. Proc. Ser. A Math. Phys. Eng. Sci.},
  454(1969):469--486, 1998.
\newblock Quantum coherence and decoherence (Santa Barbara, CA, 1996).

\bibitem{Revuz}
D.~Revuz.
\newblock {\em Markov chains}.
\newblock North-Holland Publishing Co., Amsterdam, second edition, 1984.

\bibitem{Spi}
F.\ Spitzer.
\newblock {\em Principles of random walk, 2nd edition}.
\newblock Springer-Verlag, Berlin, 2001.

\bibitem{Sym}
K.~Symanzik.
\newblock Euclidean quantum field theory.
\newblock In {\em Local quantum field theory (Varenna 1968)}. Academic Press,
  New York, 1969.

\bibitem{HofHolSla}
Remco van~der Hofstad, Frank den Hollander, and Gordon Slade.
\newblock A new inductive approach to the lace expansion for self-avoiding
  walks.
\newblock {\em Probab. Theory Related Fields}, 111(2):253--286, 1998.

\bibitem{Vermet}
Franck Vermet.
\newblock Phase transition and law of large numbers for a non-symmetric
  one-dimensional random walk with self-interactions.
\newblock {\em J. Appl. Probab.}, 35(1):55--63, 1998.

\bibitem{Westwater-I}
M.~J. Westwater.
\newblock On {E}dwards' model for long polymer chains.
\newblock {\em Comm. Math. Phys.}, 72(2):131--174, 1980.

\bibitem{WilCle}
C.~P. Williams and S.~H Clearwater.
\newblock {\em Explorations in quantum computing}.
\newblock Springer-Verlag-Telos, New York, 1998.

\bibitem{Woe}
W.\ Woess.
\newblock {\em Random walks on infinite graphs and groups}.
\newblock Cambridge University Press, Cambridge, 2000.

\end{thebibliography}
 }

\bigskip
\noindent
Massimo \textsc{Campanino}\\
Dipartimento di Matematica\\
 Universit\`a degli Studi di Bologna\\
piazza di porta San Donato 5\\
I-40126 Bologna, Italy\\
E-mail: campanin@dm.unibo.it

\bigskip
\noindent
Dimitri \textsc{Petritis}\\
Institut de Recherche Math\'ematique de Rennes\\ 
Universit\'e de Rennes I and CNRS (UMR 6625)\\
Campus de Beaulieu\\
F-35042 Rennes Cedex, France\\
E-mail: Dimitri.Petritis@univ-rennes1.fr

\end{document}